\documentclass[12pt]{article}
\usepackage{amsmath}
\usepackage[all]{xy}
\usepackage{amsfonts}
\usepackage{amssymb}
\usepackage{amscd}
\usepackage{amsthm}
\usepackage{latexsym}
\usepackage{amsbsy}

\newtheorem{lem}{Lemma}[section]
\newtheorem{prop}{Proposition}[section]
\newtheorem{theorem}{Theorem}[section]

\newtheorem{defin-prop}{Definition-Proposition}[section]
\newtheorem{defin}{Definition}[section]
\newtheorem{example}{Example}[section]

\newcommand{\ten}{\otimes_R}

\newcommand{\ot}{\otimes}

\newcommand{\sw}[1]{^{(#1)}}
\newcommand{\Hc}{\mathcal{H}}

\newcommand{\C}{\mathcal{C}}
\newcommand{\tene}{\otimes_{\re}}

\newcommand{\re}{R^{\text{e}}}
\newcommand{\en}[1]{#1\ot #1^{\text{op}}}
\newcommand{\En}[1]{#1^{\text{e}}}
\newcommand{\ra}{\rightarrow}
\newcommand{\lra}{\longrightarrow}
\newcommand{\mt}{\mapsto}

\newcommand{\D}{\mathcal{D}}
\newcommand{\cp}{\text{cop}}
\newcommand{\op}{\text{op}}
\newcommand{\bbc}{\mathbb{C}}

\newcommand{\CE}{\C ^{\text{e}}}

\begin{document}
\title
{Cyclic Cohomology of Corings}

\author{  Bahram Rangipour   \\ \\
    Department of Mathematics \\
    The Ohio State University \\
    Columbus, OH 43210, USA
} \maketitle

\begin{abstract}
We define cyclic cohomology of corings over not  necessarily
commutative algebras.  We observe that being para Hopf algebroid of
enveloping algebra of an algebra is   the key fact which allows us
to define this cohomology. This observation enables us to define
Hopf cyclic cohomology of corings on which para Hopf algebroids act.

\end{abstract}
\maketitle
\section{Introduction} Cyclic cohomology invented by  Connes \cite{c}
as  a noncommutative analogue of de Rham homology, it is a
far-reaching part of noncommutative geometry which has successfully
participated in many areas  such as Novikov conjecture for
hyperbolic groups and noncommutative index theory.   Connes defines
it as the   cohomology of the invariant co-chains under the action
of the cyclic operator which form a subcomplex of Hochschild complex
\cite{c}. Cyclic cohomology can be defined, in the most general
definition, for any cocyclic object in any abelian category
\cite{c-ext}. For a comprehensive account in cyclic (co)homology one
can see  \cite{ld}. A cocyclic object in a category is a covariant
functor from the cyclic category to that category. Cyclic cohomology
as such is a very powerful cohomology theory with all good
properties that a cohomology theory could have. Its applications and
importance   in noncommutative geometry urges us to study and
develop it to new algebraic constructions.

 Hopf cyclic  cohomology of Hopf algebras  was defined  by Connes and Moscovici
for any  Hopf algebra satisfying  some certain properties
\cite{achm98}. They showed that cyclic cohomology of their Hopf
algebra $\Hc_n$, assigned to any manifold of dimension $n$,
determines the class of the index cyclic cocycle which computes the
index of transversal elliptic operator on that manifold.  It is
shown that this cohomology is a very fit generalization of Lie
algebra cohomology and group cohomology to the realm of quantum
groups \cite{achm98, new}. Cyclic cohomology of algebras, Hopf
cyclic cohomology of Hopf algebras, and almost all other cyclic
cohomologies known so far has been unified in a single theory called
Hopf cyclic cohomology
 with coefficients \cite{hkrs}. Hopf cyclic cohomology is also seen as  a
  generalization of invariant de Rham cohomology in
noncommutative geometry \cite{inv}.

The role of groups in geometry as the source of symmetry is played
by  Hopf algebras in noncommutative geometry. The need of groupoids
in classical geometry  as a group with several objects repeats in
noncommutative geometry to introduce Hopf algebras with several
objects or Hopf algebroids.  We refer the reader to
\cite{win,e-n,hub,lu,ma,mal,ra,rin,Xu} to trace the evolution of
Hopf algebroids.  In the first attempt to develop Hopf cyclic
cohomology for Hopf algebroids, para Hopf algebroids are introduced
as the suitable gadgets on which cyclic cohomology of Hopf algebras
can be generalized  \cite{para}. A para Hopf algebroid is a
bialgebroid endowed with anti-algebra map called twisted antipode
which satisfies some certain conditions.  An example of para Hopf
algebroid is
   Connes-Moscovici  algebra in non-flat case. The other
   examples are   some  special crossed product algebras such as quantum
   torus.

    Corings can be
   seen as  coalgebras over noncommutative algebras.
   Their examples are  speared  from  the algebras of
   differential forms over a manifold  to Hopf-Galois
   extensions. One can assign to any ring extension a
   coring which speak of the fact that coring is an algebraic object,
   on the other hand to any finitely generated  projective
  module over an algebra one assigns a coring which turns
   to the coring of $n$ by $n$ matrices over the algebra in the case
   that the module is free. The latter is an example
   of crossed product coring.  Many cohomology theories such as
  intertwine cohomology can be redefine in terms of
  cohomology of coring \cite{b-w}. In this paper we develop
  Hochschild and cyclic cohomology of coring. We observe that the
  idea of Hopf-cyclic cohomology is very useful to pass the
  difficulty that the cyclic operator is not well-defined when the
  ground algebra is not commutative and/or the left and right
  action of ground algebra are not the
  same.To pass this difficulty one has to understand the real
  nature of Hochschild and cyclic complex of a coalgebra. Although
  for a coalgebra all the tensors are over the ground field,
  the role of ground field in the first tensor  is
  different from the other tensors. In the first tensor it plays as a para
  Hopf algebroid and in the other tensors its role is the same as an
  algebra.

Although the concept of para Hopf algebroid and in general
bialgebroid is not, in contrast with that of  bialgebra, self-dual,
we show that the dual theory of its cyclic cohomology is indeed rich
and contains  some algebraic information. Para Hopf algebras can be
seen as a very richer  source of symmetry in noncommutative geometry
than Hopf algebras or quantum groups, they can  act on coalgebras
and algebras with several objects, which are called coring and ring
in this paper respectively. Throughout this paper we always denote
ground fields by $k$,  ground algebras by $R$ and $S$, corings by
$\C $ and $\D$,  para Hopf algebroids by $H$, and a Hopf algebras by
$\mathcal{H}$.


\section{Corings and para Hopf algebroids}
In this section we introduce corings and their sources of symmetry
which is the category of  para Hopf algebroids. They can be seen as
coalgebras and Hopf algebras with several objects respectively. We
give some examples which will be used  in the sequel sections.
\begin{defin}\rm{
  A
coring over a $k$-algebra $R$ is an $(R,R)$-bimodule $\C$ with two
$(R,R)$-bimodule  maps
$$ \Delta: \C \rightarrow \C\ten \C, \qquad {\text{ and }}\qquad\epsilon:\C\rightarrow R, $$
that are called coproduct and counit respectively, with the
following properties
$$ (Id_\C\ten \Delta)\circ\Delta=
(\Delta\ten Id_\C)\circ\Delta,\quad \text{and}\quad (Id_\C\ten
\epsilon)\circ \Delta=Id_\C=(\epsilon\ten Id_\C)\circ\Delta.$$}

\end{defin}
We use Sweedler's notation for working with the image  of $\Delta$
in $\C\ten\C$, i.e.,  for $c\in \C$ we write  $\Delta(c)=\sum
c\sw{1}\ten c\sw{2}$, when the summation is understood. Sometimes we
even omit   the symbol  $\sum$.

One can think of corings as a generalization of coalgebras and
algebras simultaneously. If you have a $k$-coalgebra it is obvious
that it is a coring over $k$ and if you  have an algebra it is a
coring over itself with identity map as coproduct and counit. We
refer the  reader to \cite{b-w}  for a comprehensive account to have
more information about corings. A coring  over $R$ can also be
interpreted as coalgebra object in the category of $R$-module.

We would like to endow the category of coring with the most general
symmetry  that  we can. We know that for  coalgebras the most
general symmetry that one can define is the action and coaction of
bialgebras. For corings we have a thick version of bialgebras as
symmetry which are called  bialgebroid. In the following we recall
the definition of bialgebroid.

Let $H$ and $R$ be two algebras  with an algebra homomorphism
$\alpha :R \rightarrow H$, and  an algebra antihomomorphism $\beta
: R \rightarrow H$ such that the images of $\alpha $ and $\beta $
commute in $H$, i.e. for all $a$, $b$ in $R$
$$\alpha(a)\beta(b)=\beta(b)\alpha(a). $$
 It follows that $H$ has an  $R$-bimodule structure defined by
\begin{center}
$axb=\alpha(a)\beta(b)x$~~~~ $\forall a,b\in R$, ~$x\in H.$\\
\end{center}
We call $(H,R,\alpha,\beta,\Delta,\epsilon)$ a bialgebroid
 if
 \begin{itemize}
\item[i)] $(H,R,\Delta, \epsilon)$ is a coring over $R$, when $H$
is a $R$ bimodule via $\alpha$ and $\beta$.
\item[ii)]
Compatibility with the  product: for all $a$, $b\in H$ and $r\in
R$,
$$\Delta(a)(\beta(r)\otimes 1-1\otimes\alpha(r) )=0 ~~~ \text{in}~ H\otimes_RH, $$
$$\Delta (ab)=\Delta(a)\Delta(b).$$
In the first relation the natural right action of $H\otimes H$ on
$H\otimes_R H$ defined by $(a\otimes_R b)(a'\otimes
b')=aa'\otimes_R bb'$ is used. While $H\otimes_R H$ need not be an
algebra, it can be easily checked that the left annihilator of the
image of $\beta \otimes 1-1\otimes \alpha $ is an algebra. Hence,
by the first relation,  the multiplicative property of $\Delta$
makes sense.

\item[iii)] The counit is unital,$$\epsilon(1)=1.$$
\end{itemize}

\begin{defin}{\rm
 A bialgebroid $(H,R,\alpha , \beta, \Delta ,\epsilon)$
is called  a {\rm Para Hopf algebroid}
 if there is an antialgebra map $T:H\rightarrow H$, called a {\rm para-antipode},
 satisfying the following conditions:
\begin{itemize}
\item[PH1)]\label{one}$T\beta = \alpha $.
\item[PH2)]\label{three}$m_H(T\otimes id)\Delta =\beta \epsilon
T:H\rightarrow H$, where $m_H : H\otimes H \rightarrow H$ is the
multiplication map of $H$.
 \item[PH3)] $T^2=id_H$, and for all $h
\in H$
\begin{equation}\label{coalgebra}
 T(h\sw{1})\sw{1}h\sw{2}\ten
T(h\sw{1})\sw{2}=1\otimes_RT(h).
\end{equation}
\end{itemize}}
\end{defin}

 Let $\C$ be a coring over an algebra $R$. We say a bialgebroid
$(H,R)$ acts on $\C$\;  if\; $\C$ is a $H$-module, here $H$ is
thought only as an algebra over $k$. For defining  the notion of
\emph{Hopf module coring}\; we need to restrict  $H$-module
structure of $\C$ a little bit more. That is, the actions of $R$
induced by $\alpha$ and $\beta$ on $\C$ have to be coincident with
the  left and right actions of $R$ respectively. By the above
assumption we can define the diagonal action of $(H,R)$ on $(\C,R)$
as follows
\begin{center}
$H\ot \C\ten \C\longrightarrow \C\ten \C$\\
$h\cdot(c\ten d)=h\sw{1}c\ten h\sw{2}d$\quad for any $h\in H$ and
$c,d\in\C$.
\end{center}
This defines an action because of   (ii) in  the definition of
bialgebroids. We know that the base algebra as such is a coring over
itself. We want to endow it an structure of $H$-module coring. To
this end,   we have to assume one more axiom for bialgebroids, that
is for any $g,h\in H$
$$\epsilon(hg)=\epsilon(h\beta(\epsilon(g))).$$
This condition in many sources is one of the axioms of
bialgebroids.  One can show that for a para Hopf algebroid this
condition is satisfied for free \cite{para}.  Thanks to the above
assumption on $\epsilon$ one  defines the following action of a
bialgebra $(H,R)$ on $R$.
\begin{center}
$H\ot R\longrightarrow R$\\
$h\cdot r=\epsilon(h\alpha(r))=\epsilon(h\beta(r))$
\end{center}
Which is a generalization of trivial action of Hopf algebras on
the ground ring.
\begin{defin}{\rm
Assume that  $(H,R)$ acts on a coring $\C$ and this action has the
above property, i.e. the left and right actions of $R$ coincide with
its actions induced by the action of $H$ on $C$. We call $\C$ an
$H$-module coring if $\Delta:\C\rightarrow \C\ten \C$ and $\epsilon:
\C\rightarrow R$ are $H$-module maps, where $\C\ten \C$ and $R$ are
$H$-module as above.}
\end{defin}

The above definition generalizes  Hopf module coalgebras. One sees
that  action of any bialgebroid on itself by multiplication
satisfies the above definition.
\begin{example} {\rm Let  $B\rightarrow A$ be  a $k$ algebra extension.
 The Sweedler's bialgebroid $A_B^{\text{e}}=A\ot_B A^{\text{op}}$ is defined as follows.
\begin{align*}
&\alpha: A \ra A_B^{\text{e}}, &&{\rm{by}} &&& a\mt a\ot_B 1_A\\
&\beta: A \ra A_B^{\text{e}}, &&{\rm{by}} &&& x\mt 1_A\ot x\\
&\Delta :A_B^{\text{e}}\lra A_B^{\text{e}}\ot_AA_B^{\text{e}} &&
{\rm by} &&& a\ot x\mt a\ot 1_A\ot_A
1_A\ot x\\
&\epsilon: A_B^{\text{e}}\lra A && {\rm by}&&& a\ot x\mt ax
\end{align*}
The Sweedler's bialgebroid $\re=R_k^{\text{e}}$ acts on any coring
$\C$ over the algebra $R$ by $$(r\ot s)c=r\cdot c\cdot s.$$ One
can show that in fact $\C$ is a $\re$-module coring.}
\end{example}

\begin{example}
Any bialgebroid $(H,R)$ acts on $A^{\text{e}}$ via $h(a\ot
b)=\epsilon_H(h\sw{1}\alpha(a))\ot \epsilon_H(h\sw{1}\beta(a))$
and makes $A^{\text{e}}$ a left $H$-module coring.
\end{example}

Let $\mathcal{H}$ be a Hopf algebra with an invertible antipode. Let
$M$ be a left $\mathcal{H}$-module and a right $\mathcal{H}$
comodule.  It is called  Yetter-Drinfeld module if
$$\Delta_M(hm)= h\sw{2}m\sw{0}\ot h\sw{3}m\sw{1}S^{-1}(h\sw{1}).$$
Then $_\mathcal{H}\!\mathcal{YD}^\mathcal{H}$, the category of
 Yetter-Drinfeld module over $\mathcal{H}$,  consists of the above
objects and usual morphisms. It  forms a braided category with the
following braiding map:
$$\sigma:M\ot N\rightarrow N\ot M,\hspace{1cm} (m\ot n)\mapsto n\sw{0}\ot n\sw{1}m.$$
Let  $A\in\; _\mathcal{H}\!\mathcal{YD}^\mathcal{H}$ be a {\it
braided commutative} algebra in
$_\mathcal{H}\!\mathcal{YD}^\mathcal{H}$, i.e., it is a left
$\mathcal{H}$ module algebra and right $\mathcal{H}^{\text{op}}$
comodule algebra  and for any $a,b\in A$ we have
$$b\sw{0}(b\sw{1}a)=ab$$

It is shown in \cite{bm} that the crossed product algebra $H=A\sharp
\mathcal{H}$ is a Hopf algebroid over $A$ with the following
structure.
\begin{align*}
&\Delta :H\rightarrow H\ot_A H, && \Delta(a\#h)=a\#h\sw{1}\ot_A
1\#
h\sw{2},\\
&\epsilon:H\rightarrow A, &&\epsilon(a\# h)=a\epsilon(h),\\
&\alpha: A\rightarrow H, &&\alpha(a)=a\# 1,\\
& \beta: A\rightarrow H,&&\beta(b)=b\sw{0}\# b\sw{1},\\
& \tau: H\rightarrow H, && \tau(a\#
h)=S(h\sw{2})S^2(a\sw{1})a\sw{0}\# S(h\sw{1})S^2(a\sw{2}).
\end{align*}

Let $S^2=id$ and $A$ is stable, i.e., $a\sw{1}a\sw{0}=a$.  We show
that the crossed product algebra $H=A\# \mathcal{H}$ is a coring
over $A$. We keep all bialgebroid structure of $H$ over $A$
defined in the above and simplify  its antipode to $T:H\rightarrow
H$ defined by
 $$T(a\# h)= S(h\sw{2})a\sw{0}\# {S}(h\sw{1})a\sw{1}$$
 \begin{prop}
 With the above definition $(H,A,\alpha,\beta,\Delta,\epsilon,T)$
 is a para Hopf algebroid.
 \end{prop}

\begin{proof}
We know it is a  bialgebroid. We show that $T^2=id$,
 and $T$ satisfies the condition (\ref{coalgebra}). To show the former,
since $T$ is an antialgebra  map and for any $h\in \mathcal{H}$ we
have $T^2(1\# h)=1\# S^2(h)$, it suffices to show $T^2(a\# 1)=a\#
1$ for any $a\in A$.
\begin{align*}
&T^2(a\# 1)=T(a\sw{0}\# a\sw{1})=T(1\# a\sw{1})T(a\# 1)=(1\#
S(a\sw{2}))(a\sw{0}\#a\sw{1})\\
&=S(a\sw{3})a\sw{0}\# S(a\sw{2})a\sw{1}=S(a\sw{1})a\sw{0}\#1=a\#1
\end{align*}
Now let us check the condition (\ref{coalgebra}). It is shown that
the condition (\ref{coalgebra}) is multiplicative \cite{para}. So we
may check it for $1\# h$ and $a\# 1$ separately. We just check the
latter and leave the other to the reader.
\begin{align*}
&T((a\#1)\sw{1})\sw{1}(a\#1)\sw{2}\ot_A
T((a\#1)\sw{1})\sw{2}\\
&=(a\sw{0}\# a\sw{1})\sw{1}(1\# 1)\ot_A(a\sw{0}\#
1)\sw{2}=(a\sw{0}\# a\sw{1})\ot_A 1\# a\sw{2}\\
&=\beta(a\sw{0})(1\#1)\ot_A 1\# a\sw{1}=1\#1\ot_A
 a\sw{0}\#a\sw{1}=1\#1\ot_A T(a\# 1)
\end{align*}
\end{proof}
\begin{example}
Let $\mathcal{H}$ be a Hopf algebra with $S^2=id$. Let
$A=\mathcal{H}$ with the conjugation action, i.e., $ha=
h\sw{1}aS(h\sw{2})$. The right coaction that we propose is
$\Delta_A(a)= a\sw{2}\ot S(a\sw{1})$. One can check that $A$ is a
braided commutative algebra in
$_\mathcal{H}\!\mathcal{YD}^\mathcal{H}$. So $A\# H$ is a para
Hopf algebroid on $A$.
\end{example}
\begin{example}
Let a discrete group $\Gamma$  acts from right by diffeomorphisms
on a manifold $M$. One can check that $A=C^\infty(M)$ with trivial
coaction of $\mathcal{H}=\bbc\Gamma$ satisfies the condition of
the above proposition, hence $C^\infty(M)\rtimes \Gamma$ is a para
Hopf algebroid over $C^\infty(M)$.
\end{example}
Let $\C$ and $\D$ be  corings over $R$  and $S$ respectively. It
is easy to check that $\C\ot \D$ is a coring over $R\ot S$, where
$\C\ot \D$ is $R\ot S$ bimodule via $(r\ot s)(c\ot d)=rc\ot sd$,
and $(c\ot d)(r\ot s)=cr\ot ds$, $c\in\C$, $d\in \D$, $r\in R$,
and  $s\in S$. Its coproduct and counit are defined by

\begin{align*}
&\Delta: \C\ot \D\rightarrow \C\ot \D\ot_{R\ot S} \C\ot \D,
&&\Delta(c\ot d)= c\sw{1}\ot d\sw{1}\ot_{R\ot S} c\sw{2}\ot
d\sw{2}\\
&\epsilon: \C\ot \D\rightarrow R\ot S, && \epsilon(c\ot
d)=\epsilon(c)\ot\epsilon(d).
\end{align*}
On the other hand  $\C^\cp$ which is $\C$ as a vector space has  a
coring structure over $R^\op$ as follows. The $R^\op$ bimodule
structure of $\C^\cp$ is defined in the usual way, i.e., $r_1
 cr_2:=r_2cr_1$.The coproduct and counit is defined by
 \begin{align*}
 &\Delta: \C^\cp\rightarrow \C^\cp\ot_{R^\op} \C^\cp,
 &&\Delta_{\C^\cp}(c)={c\sw{2}}\ot_{R^\op} {c\sw{1}}\\
 &\epsilon_{\C^\cp}:\C^\cp\rightarrow R^\op, &&
 \epsilon_{\C^\cp}=\epsilon(c).
 \end{align*}
Now we can define enveloping coring  of a coring $(\C,R)$, which
is $(\C\ot \C^\cp, R\ot R^\op)$ defined as above. We denote this
coring by $(\CE, R^{\text{e}})$.
\section{Hopf cyclic co/homology of corings}
In this section we recall the cyclic cohomology of coalgebras  and
then we will extend this theory to the category of corings. We will
observe that the key concept that enables us to proceed  is the
notion of para Hopf algebroid.

 Let $C$ be a coalgebra over a commutative ring $k$. Recall that
 the cocyclic module assigned to $C$ is
 $\{C^n=C^n(C)=C^{\ot n+1}\}_{n\ge 0}$ with the following
 maps:\\[.2cm]
$\delta_i:C^n\rightarrow C^{n+1},\;\;\; i=0,1,\dots,n+1$\\
$\delta_i(c_0\ot c_1\ot\dots\ot c_n)=(c_0\ot c_1\ot \dots\ot\Delta(c_i)\ot\dots\ot c_n)$\\
$\delta_{n+1}(c_0\ot c_1\ot \dots\ot c_n)=(c_0\sw{2}\ot c_1\ot\dots\ot c_n\ot c_0\sw{1})$\\[.5cm]
$\sigma_j:C^n\rightarrow C^{n-1}\;\; j=0,1,\dots,n-1$\\
$\sigma_j(c_0\ot c_1\ot\dots\ot c_n)=(c_0\ot c_1\ot\dots\ot\epsilon(c_{i+1})\dots\ot c_n)$\\[.5cm]
$\tau:C^n\rightarrow C^n$\\
$\tau(c_0\ot c_1\ot\dots\ot c_n)=(c_1\ot\dots\ot c_n\ot c_0)$

Now let us try to extend the above structure for a  coring $(\C,R)$.
If we simply change the tensor over $k$ to tensor over $R$, the
first problem that we face is that the cyclic operator $\tau$ and
the last face $\delta_{n+1}$ are not well-defined. To fix this
problem we observe that if we change the first tensor  to tensor
over $R^e$ and the others to tensor over $R$ then the operators are
well-defined and we have the following proposition.

\begin{prop}
Let $(\C,R)$ be a coring then the following define a cocyclic
$k$-module.

$\{C^n=C^n(\C)=\C\tene\C^{\ten n}\}_{n\ge 1}$, $C^0=\C\tene R$
with the following
 maps:\\[.2cm]
$\delta_i:C^n\rightarrow C^{n+1},\;\;\; i=0,1,\dots,n+1$\\
$\delta_i(c_0\tene c_1\ten\dots\ten c_n)=(c_0\tene c_1\ten \dots\ten\Delta(c_i)\ten\dots\ten c_n)$\\
$\delta_{n+1}(c_0\tene c_1\ten \dots\ten c_n)=(c_0\sw{2}\tene c_1\ten\dots\ten c_n\ten c_0\sw{1})$\\[.5cm]
$\sigma_j:C^n\rightarrow C^{n-1}\;\; j=0,1,\dots,n-1$\\
$\sigma_j(c_0\tene c_1\ten\dots\ten c_n)=(c_0\tene c_1\ten\dots\ten\epsilon(c_{i+1})\dots\ten c_n)$\\[.5cm]
$\tau:C^n\rightarrow C^n$\\
$\tau(c_0\tene c_1\ten\dots\ten c_n)=(c_1\tene\dots\ten c_n\ten
c_0)$
\end{prop}
\begin{proof} See the proof of Theorem \ref{th} as a more general
case. \end{proof}

 We denote its Hochschild, cyclic and periodic
cyclic cohomology by $HH^\ast(\C)$, $HC^\ast(\C)$, and
$HP^\ast(\C)$ respectively.

\begin{lem}
Let $(\C,R)$ be a $H$-coring module. Then for any $r,s\in R$, and
$c\in \C$ we have

$$\alpha(r)c\ot_H s=c\ot_Hsr,\;\; \text{and}\;\; \beta(r)c\ot_H
s=c\ot_Hrs.$$

\end{lem}

\begin{proof}
Using the fact that $T^2=id_H$ and $T\beta=\alpha$, we conclude
that $T\alpha=\beta$. So $\alpha(r)c\ot_H s=c\cdot\beta(r)\ot_H
s=c\ot_H\beta(r)\cdot
s=c\ot_H\epsilon_H(\beta(r)\beta(s))=c\ot_H\epsilon_H(\beta(sr))=c\ot_Hsr.$
And similarly,  $\beta(r)c\ot_H s=c\cdot\alpha(r)\ot_H
s=c\ot_H\alpha(r)\cdot
s=c\ot_H\epsilon_H(\alpha(r)\alpha(s))=c\ot_H\epsilon_H(\alpha(rs))=c\ot_Hrs.$
\end{proof}
\begin{lem}\label{action}
For any $n\ge 1$, $h\in H$, and $c_1,\dots ,c_n\in \C$ one has:
\begin{align*}
&T(h\sw{1})\sw{1}h\sw{2}c_1\ten\dots\ten
T(h\sw{1})\sw{n-1}h\sw{n}c_{n-1}\ten
T(h\sw{1})\sw{n}c_n=\\
&c_1\ten\dots\ten c_{n-1}\ten T(h)c_n
\end{align*}
\begin{proof}
First, we see that the following defines a well-defined map,
\\[.1cm]
$\Psi: H^{\ten n}\ot\C^{\ot n}\rightarrow \C^{\ten n} $\\[.1cm]
$\Psi(h_1\ten\dots\ten h_n\ot c_1\ot\dots\ot
 c_n)=h_1c_1\ten\dots\ten h_nc_n$.\\[.1cm]
 Next, we apply $\Delta_H$, the comultiplication  of $H$, $n-2$ times on the two sides of the
 crucial equation \eqref{coalgebra} to get
 $$T(h\sw{1})\sw{1}h\sw{2}\ten\dots\ten
T(h\sw{1})\sw{n-1}h\sw{n}\ten T(h\sw{1})\sw{n}=1\ten\dots\ten 1\ten
T(h).$$ Then we apply both hand sides on $c_1\ten\dots\ten c_n$.

\end{proof}

\end{lem}
The reader may ask what algebra except  $\re$ can be used to
balance the above complex to get a cocyclic module. This question
is answered in the following theorem.

In the following theorem \;$\C$\;  is  considered as a right $H$
module via  the antipode of $H$ and $\C^{\ten n}$ is left $H$ module
diagonally.
\begin{theorem}\label{th}
Let $(H,R)$ be a para Hopf algebroid  and $(\C,R)$ an $H$-module
coring. Then the following define a cocyclic module:

$\{C^n=C^n(\C)=\C\ot_H\C^{\ten n}\}_{n\ge 0}$ with the following
 maps:\\[.2cm]
 $\delta_0,\delta_1:C^0\ra C^1$,  defined by, \\
 \begin{align*}
 &\delta_0(c\ot_H r)=c\sw{1}\ot_H
 \beta(r)c\sw{2}\\
 &\delta_1(c\ot_H r)=c\sw{2}\ot_H\alpha(r)c\sw{2}
 \end{align*}
$\delta_i:C^n\rightarrow C^{n+1},\;\;\; i=0,1,\dots,n+1$
\begin{align*}
&\delta_i(c_0\ot_H c_1\ten\dots\ten c_n)=(c_0\ot_H c_1\ten \dots\ten\Delta(c_i)\ten\dots\ten c_n)\\
&\delta_{n+1}(c_0\ot_H c_1\ten \dots\ten c_n)=(c_0\sw{2}\ot_H
c_1\ten\dots\ten c_n\ten c_0\sw{1})
\end{align*}
$\sigma_j:C^n\rightarrow C^{n-1}\;\; j=0,1,\dots,n-1$
\begin{equation*}
\sigma_j(c_0\ot_H c_1\ten\dots\ten c_n)=(c_0\ot_H
c_1\ten\dots\ten\epsilon(c_{i+1})\ten\dots\ten c_n)
\end{equation*}
$\tau:C^n\rightarrow C^n$
\begin{equation*}
\tau(c_0\ot_H c_1\ten\dots\ten c_n)=(c_1\ot_H\dots\ten c_n\ten c_0)
\end{equation*}
\end{theorem}

\begin{proof}
We should show that $\delta_i, \sigma_j$, and $\tau$ satisfy the
cocyclic module relations.

 First cosimplicial module relations:
\begin{center}
$\delta_i\delta_j=\delta_{j}\delta_{i-1}$ $j< i$,\\
$\sigma_i\sigma_j=\sigma_j\sigma_{i+1}$ $j\le i$,\\
$\sigma_j\delta_i=\left\{\begin{matrix}
\delta_i\sigma_{j-1}\;\;&\text{if}\;\;i<j\\ id
\;\;&\text{if}\;\;i=j
\;\;\text{or}\;\; i=j+1\\

\delta_{i-1} \sigma_j\;\;&\text{if}\;i>j+1,\end{matrix}\right.$

\end{center}
and relations of cyclic operator with faces and degeneracies:
$\tau\delta_i=\delta_{i-1}\tau$,  $\tau\sigma_j=\sigma_{j-1}\tau$
and the crucial relation $\tau^{n+1}=id$.

We leave to the reader to check the the above relations because
they are the same as the cyclic structure of coalgebras but we
show the above maps are well-defined. It is easy to see that
$\delta_i$, $i=0,\dots,n$ and $\sigma_j$,  $j=0,1,\dots,n-1$ are
well-defined due to the fact that $\Delta_\C$  and $\epsilon_\C$
are $ H$-linear.
 If we prove $\tau$ is well-defined, then since
 $\delta_{n+1}=\tau\delta_0$ we have done this stage completely.

In degree zero the cyclic operator is $id_{C^0}$ and hence  it is
well-defined.  Let $n\ge 1$,  $h\in H$, $r\in R$ and
$c_0,\dots,c_n\in \C$.\\  Using Lemma \ref{action} we have
\begin{align*}
&\tau(c_0\ot h\cdot c_1\ot \dots\ot
c_n)=\tau(c_0\ot h\sw{1}c_1\ot\dots\ot h\sw{n}c_n)\\
&=h\sw{1}c_1\ot_Hh\sw{2}c_2\ten\dots\ten h\sw{n}c_n\ten c_0\\
&=c_1T(h\sw{1})\ot_Hh\sw{2}c_2\ten\dots\ten h\sw{n}c_n\ten c_0\\
&=c_1\ot_H T(h\sw{1})\sw{1}h\sw{2}c_2\ten\dots\ten
T(h\sw{1})\sw{n-1}h\sw{n}c_n\ten T(h\sw{1})\sw{n}c_0\\
&=c_1\ot_H c_2\ten\dots\ten c_n\ten T(h)c_0=\tau(c_0h\ot c_1\dots\ot
c_n).
\end{align*}

There are also some ambiguities over $\ten$'s which we should get
rid of. We show $\tau$ is well-defined over the first $\ten$ and
leave the rest to the reader.
 \begin{align*}
 &\tau(c_0\ot  c_1\cdot
r\ot\dots\ot c_n)=\\
&\beta(r)c_1\ot_H c_2\ten\dots\ten c_n\ten c_0=\\
&c_1\cdot\alpha( r)\ot_H c_2\ten\dots\ten c_n\ten c_0=\\
&c_1\ot_H \alpha(r)\triangleright c_2\ten\dots\ten c_n\ten
c_0=\\
&c_1\ot_H \alpha(r)c_2\ten\dots\ten c_n\ten c_0=\\
&\tau(c_0\ot c_1\ot r\cdot c_2\ot\dots\ot c_n).
\end{align*}
\end{proof}

The Hopf  cyclic cohomology of  $(\C, R)$ under the action of $H$ is
denoted by $HC_H^*(\C, R)$.

\begin{defin}{\rm
We say a $H$-module  coring $(\C, R)$ is $H$ coseparable if there is
a $H$-linear map $\delta: \C\ten \C \rightarrow R$ satisfying
\begin{enumerate}
    \item $\delta\Delta=\epsilon $
    \item $(id_\C\ten\delta)(\Delta\ten id_\C)=(\delta\ten
    id_\C)(id_\C\ten\Delta)$
\end{enumerate}}
\end{defin}

 For example if  an algebra extension $B\hookrightarrow A$ is split, i.e.,
there is a $B-B$ bimodule map $E:A\rightarrow B$ such that $E(1)=1$,
then the Sweedler coring $A\ot_BA$ is coseparable \cite{guz}.

\begin{defin}{\rm (Haar system for bialgebroids) Let $(H,R)$ be a bialgebroid.
Let $\theta : H\longrightarrow R$ be a  right R-module map. We
call $\theta$ a left
 Haar system for $H$ if  for all $h\in H$
$$\sum\alpha(\theta(h^{(1)})) h^{(2)}=\beta(\theta(h))1_H.$$
We call $\theta$ a normal left Haar system if $\theta(1_H)=1_R.$}
\end{defin}
\begin{prop}
If a para Hopf algebroid $H$ is coseparable as a module coring over
itself, then it admits a normalized  Haar system.
\end{prop}

\begin{proof}

 Let $\delta: H\ten H\rightarrow R$ be a
coseparating map for $H$. We define $\theta: H\rightarrow R $ by
$\theta(h)=\delta(1,h)$. We see that $\theta(h\sw{1})\ten
h\sw{2}=\delta(1,h\sw{1})\ten h\sw{2}=1\ten \delta(1\ten
h)=1\ten\theta(h).$ So $\theta$ is a Haar system. It is easily seen
that normality comes from the first condition of coseparability.
\end{proof}

\begin{prop}{\rm (\cite{para})}
Let  $(H,R)$ be a para Hopf algebroid  that admits a normal left
Haar system. Then $HC^{2i+1}_H(H, R)=0$ and $HC^{2i}_H(H,
R)=\ker(\alpha-\beta)$ for all $i\geq 0$.
\end{prop}

\begin{prop}
If $(\C,R)$ is  $H$ coseparable, then $$HC^n_H(\C,R)=\left\{%
\begin{array}{ll}
    HC_H^0(\C,R) & \hbox{n\;\;\text{even} ;} \\
    0 & \hbox{n \;\;\text{odd}.} \\
\end{array}%
\right. $$
\end{prop}
\begin{proof}
It is easily checked that the following is a  homotopy between
identity and zero map on the Hochschild complex. The rest is a
standard application of Connes $SIB$ long exact sequence.
 \begin{align*}
& h: C^n\rightarrow C^{n-1}\\
&h(c_0\ot_H c_1\ten \dots\ten c_n)=c_0\sw{1}\ot_H
\delta(c_0\sw{2}\ten c_1)\ten c_2\ten\dots\ten c_n
\end{align*}
\end{proof}

As we know corings generalize coalgebras and algebras
simultaneously, this interestingly happens for  their cyclic
co/homology too. Obviously this generalization can be seen for
coalgebras but for algebras we need first to see how the dual of the
above cocyclic module looks like. We recall from \cite{c}
 see also \cite{dual}, roughly speaking, dual of a cocyclic module, say $C^n$, is a cyclic
 module with the same vector space structure
whose faces are degeneracies of $C^n$ plus the extra degeneracy
which is defined by $\sigma_{-1}=\tau\sigma_{n-1}$. The
degeneracies of the dual cyclic module are faces of the original
cocyclic module, and finally    the cyclic operator of the dual
cyclic module is the inverse of the cyclic operator of the
original cyclic module. In the following we denote the dual of the
above cocyclic module by $\hat C_n(\C)$. Its cyclic  structure is
given by:

 $d_0,d_1:\hat C_1\ra \hat C_0$ defined by\\[.1cm]
 $d_0(c_0\ot_Hc_1)=c_0\ot_H\epsilon(c_1), d_1(c_0\ot_H
 c_1)=c_1\ot_H\epsilon(c_0)$\\
$d_i:\hat C_n\rightarrow \hat C_{n-1},\;\;\; i=0,1,\dots,n$\\
$d_i(c_0\ot_H c_1\ten\dots\ten c_n)=(c_0\ot_H c_1\ten
\dots\ten\epsilon(c_{i+1})\ten\dots\ten
c_n)$\\[.1cm]
$d_n(c_0\ot_H c_1\ten\dots\ten c_n)=(c_n\ot_H \epsilon(c_0)\ten \dots\ten c_{n-1})$\\[.5cm]
$s_j:\hat \C_n\rightarrow \hat C _{n+1}\;\; j=0,1,\dots,n$\\[.1cm]
$s_n(c_0\ot_H c_1\ten\dots\ten c_n)=(c_0\ot_H
c_1\ten\dots\ten\Delta(c_{i+1})\dots\ten c_n)$\\
 $s_j(c_0\ot_H
c_1\ten\dots\ten c_n)=(c_0\sw{2}\ot_H
c_1\ten\dots\ten\epsilon(c_{i+1})\dots\ten c_n\ten c_0\sw{1})$\\[.5cm]
$t:\hat C_n\rightarrow \hat C_n$\\[.1cm]
$t(c_0\ot_H c_1\ten\dots\ten c_n)=(c_n\ot_H c_0 \ten c_1\ten
\dots\ten c_n)$

\begin{prop}
Cyclic module of an algebra $A$ is isomorphic to dual Hopf-cyclic
module of the Sweedler coring $\en{A}$ as $\en{A}$- module coring.
\end{prop}
\begin{proof}
Let $C_n(A)$ denote the ordinary cyclic complex of the algebra
$A$. We define the following maps and leave to the reader to check
that these map are cyclic maps, i.e. they commute with cyclic
structure, and are inverse of each other.

$\Psi: C_n(A)\longrightarrow  \hat C_n(\En{A})$\\
$\Psi(a_0\ot\dots\ot a_n)=1_{\En{A}}\ot_{\En{A}} a_0\ot
1_A\ten\dots \ten a_{n-1}\ot 1_A\ten a_{n-1}\ot a_n $\\[.2cm]

$\Phi : \hat C_n(\En{A})\longrightarrow C_n(A)$\\
$\Phi(a_0\ot b_0\ot_{\En{A}} a_1\ot b_1\ten \dots \ten a_n\ot
b_n)=b_0a_1\ot b_1a_2\ot \dots \ot b_{n-1}a_n\ot b_na_0$

\end{proof}

Any para Hopf algebroid is  a module coring on itself via
multiplication. We show that cocyclic module of a para Hopf
algebroid defined in \cite{para} is isomorphic  to its  cocyclic
module as a module coring on itself. To this end, we first recall
cocyclic module of a para Hopf algebroid from \cite{para}.

Let $(H,R)$ be  para Hopf algebroid. Its cocyclic module as a para
Hopf algebroid, also called Connes-Moscovici cocyclic module of $H$
is defined as follows:
$$H_\natural^0=R, \text{ and} \;H_\natural^n=H\otimes_R H\otimes_R
  \dots \otimes_RH \qquad (n\text{-fold tensor product}).$$
   The cofaces $\delta_i$ and  codegeneracies  $\sigma_i$  are defined by:
\begin{eqnarray*}
\delta_0(a)=\alpha(a),~ \delta_1(a)=\beta(a)&& \text{for all}~ a\in R=H^0_\natural\\
\delta_0(h_1\otimes_R\dots \otimes_Rh_n)&=&1_H\otimes_R h_1\otimes_R\dots \otimes_Rh_n \\
\delta_i(h_1\otimes_R\dots \otimes_Rh_n)&=&
h_1\otimes_R\dots\otimes_R\Delta (h_i)\otimes_R\dots \otimes_Rh_n
  \;\;\text{for}\;\;1\leq i\leq n \\
\delta_{n+1}(h_1\otimes_R\dots \otimes_Rh_n)&=&h_1\otimes_R\dots \otimes_Rh_m\otimes_R 1_H \\
\sigma_i(h_1\otimes_R\dots \otimes_Rh_n)&=&
h_1\otimes_R\dots\otimes_R\epsilon(h_{i+1})\otimes_R\dots
\otimes_Rh_n
 \;\;\text{for}\;\;0\leq i\leq n. \\
 \end{eqnarray*}
 The cyclic operator $\tau$  defined by
$$\tau_n(h_1\otimes_R\dots \otimes_Rh_n
)=T(h_1)\vartriangleright(h_2\ten\dots \ten h_n\ten 1_H),$$

\begin{prop}
Let $(H,R)$ be a para Hopf algebroid. Then its Hopf-cyclic module
as module coring over itself is isomorphic to its Connes-Moscovici
cocyclic module.
\end{prop}
\begin{proof}
We define the following maps and leave to the reader to verify
that they are cyclic maps and inverse to one another.

\noindent $\Psi: C^n(H)\longrightarrow H_\natural^n$\\
 $\Psi (h_0\ot_H h_1\ten \dots\ten h_n)=T(h_0)\vartriangleright(h_1\ten\dots \ten
 h_n)$\\[.2cm]
 $\Phi: H_\natural^n \longrightarrow C^n(H)$\\
$ \Phi(h_1\ten \dots\ten h_n)=1_H\ot_H h_1\ten\dots\ten h_n.$
\end{proof}


\end{document}